\documentclass[11pt]{article}
\usepackage{amsfonts}
\usepackage{amssymb}
\usepackage{graphicx}
\usepackage{epsfig}
\usepackage{amsmath}

\newcommand{\ignore}[1]{}
\ignore{The text in between these parentheses is ignored by LaTeX
ignore}

\def\@begintheorem#1#2{\par\bgroup{\sc #1\ #2. }\it\ignorespaces}
\def\@opargbegintheorem#1#2#3{\par\bgroup{\sc #1\ #2\ (#3). } \it\ignorespaces}
\def\@endtheorem{\egroup}
\newtheorem{theorem}{Theorem}[section]
\newtheorem{corollary}[theorem]{Corollary}
\newtheorem{lemma}[theorem]{Lemma}
\newtheorem{proposition}[theorem]{Proposition}
\newtheorem{problem}[theorem]{Problem}
\newtheorem{example}[theorem]{Example}
\newtheorem{algorithm}[theorem]{Algorithm}
\newtheorem{definition}[theorem]{Definition}
\newcommand{\bt}[1]{\begin{theorem}\label{#1}}
\newcommand{\bc}[1]{\begin{corollary}\label{#1}}
\newcommand{\bl}[1]{\begin{lemma}\label{#1}}
\newcommand{\bp}[1]{\begin{proposition}\label{#1}}
\newcommand{\bpro}[1]{\begin{problem}\label{#1}}
\newcommand{\be}[1]{\begin{example}\rm\label{#1}}
\newcommand{\ba}[1]{\begin{algorithm}\rm\label{#1}}
\newcommand{\bd}[1]{\begin{definition}\rm\label{#1}}

\newcommand{\et}{\end{theorem}}
\newcommand{\ec}{\end{corollary}}
\newcommand{\el}{\end{lemma}}
\newcommand{\ep}{\end{proposition}}
\newcommand{\epro}{\end{problem}}
\newcommand{\ee}{\end{example}}
\newcommand{\ea}{\end{algorithm}}
\newcommand{\ed}{\end{definition}}

\def\N{\mathbb{N}}
\def\R{\mathbb{R}}
\def\Z{\mathbb{Z}}
\def \conv {\,\mbox{conv\,}}
\def \vert {\mbox{vert}}
\def \sign {\mbox{sign}}
\def \F {{{\cal F}}}

\newcommand{\boproof}{\noindent {\em Proof. }}
\newcommand{\eoproof}{\hspace*{\fill} $\square$ \vspace{5pt}}

\begin{document}

\title{\bf Nonlinear Bipartite Matching}
\author{
Yael Berstein
\thanks{
Supported in part by a Technion Graduate School Fellowship.}
\and
Shmuel Onn
\thanks{Supported in part by a grant from ISF - the
Israel Science Foundation, by the Technion President Fund,
and by the Jewish Communities of Germany Research Fund.}}

\date{}
\maketitle


\begin{abstract}
We study the problem of optimizing nonlinear objective functions
over bipartite matchings. While the problem is generally intractable,
we provide several efficient algorithms for it, including a
deterministic algorithm for maximizing convex objectives,
approximative algorithms for norm minimization and maximization,
and a randomized algorithm for optimizing arbitrary objectives.
\end{abstract}

\section{Introduction}
\label{Introduction}

Let $N:=\{(i,j):1\leq i,j\leq n\}$ be the set of edges of the complete
bipartite graph $K_{n,n}$. In this article we consider the following
broad generalization of the standard linear bipartite matching problem.

\vskip.2cm\noindent
{\bf Nonlinear Bipartite Matching}. Given positive integers $d,n$,
integer weight functions $w^1,\dots,w^d$ on $N$, and an arbitrary function
$f:\R^d\longrightarrow\R$, find a perfect matching $M\subset N$ maximizing
(or minimizing) the objective function $f(w^1(M),\dots,w^d(M))$
where $w^k(M):=\sum\{w^k(i,j):(i,j)\in M\}$.

\vskip.2cm\noindent
Identifying perfect matchings in $K_{n,n}$ with permutation matrices
and weight functions with integer matrices in the usual way, the problem
has the following nonlinear integer programming formulation:
$$\max\, \{f(w^1 x,\dots, w^d x):\ \sum_{i=1}^n x_{i,j}=1,
\ \sum_{j=1}^n x_{i,j}=1,\ x\in\N^{n\times n}\}$$
where $w^k x:=\sum_{i=1}^n\sum_{j=1}^n w^k_{i,j} x_{i,j}$ for
$k=1,\dots,d$, and where $\N$ stands for the nonnegative integers.

\noindent
The problem can be interpreted as {\em multiobjective} bipartite matching:
given $d$ different linear objective functions $w^1,\dots,w^d$, the goal is
to maximize (or minimize) their ``balancing" given by $f(w^1 x,\dots,w^d x)$.
The standard linear bipartite matching problem is the special
case of $d=1$ and $f$ the identity on $\R$.

Beyond the intrinsic interest in studying the above natural nonlinear
extension of the standard matching problem, it is interesting to consider
it in connection with its various variants and relatives in the literature,
including in \cite{CKM,DO,HOR,Kar,MVV,Onn,OR,OS,Pap,PY,YMS} and
references therein. These variants will be discussed in detail in Section 2.
In particular, in \cite{OR}, the problem of maximizing an objective
function $f(w^1(\cdot),\dots, w^d(\cdot))$ with $f$ convex and $d$ fixed
was considered for combinatorial optimization families in general. It was
shown, extending earlier results of \cite{Onn,OS}, that if the polyhedra
underlying the problem have nice edge symmetry (few edge-directions)
then the problem can be solved in strongly polynomial time. This resulted in
polynomial time algorithms for maximizing a convex $f$ for various problems
including vector partitioning, matroids, and transportation problems with
fixed number of suppliers. However, the Birkhoff polytope which underlies
our bipartite matching problem has exponentially many edge-directions
(Proposition \ref{EdgeDirections}, Section 2)
and hence the methods of \cite{Onn,OR,OS} fail.

Nonlinear bipartite matching is generally intractable, since already for
fixed $d=1$ (single weight function), the problem of minimizing a family
of very simple convex univariate functions $f_u:\R\longrightarrow\R$
defined by $f_u(y):=(y-u)^2$ with $u$ an integer parameter is NP-hard
(Proposition \ref{NPHard} part 1, Section 2). \break
Therefore, for the most part, the complexity of our results will depend
on the unary size of the weights, that is, on $\max|w^k_{i,j}|$.
In particular, our algorithms will have polynomial complexity for
{\em binary} weights, that is, with $w^k_{i,j}\in\{0,1\}$ for all $i,j,k$.
In this case, letting $E_k$ be the support of $w^k$ for each $k$, the
problem becomes that of finding a perfect matching $M\subset N$
maximizing (or minimizing) $f(|M\cap E_1|,\dots,|M\cap E_d|)$.
The problem with binary weights is not easy either: the complexity with $f$
an arbitrary function is unknown for any fixed $d\geq 2$; and for
variable $d$ it is again NP-hard for minimizing the convex multivariate
extension of $f_u$ above, i.e. the family of functions
$f_u:\R^d\longrightarrow\R$ defined by $f_u(y):=\sum_{k=1}^d(y_k-u_k)^2$
and parameterized by $u\in\Z^d$ (Proposition \ref{NPHard} part 2, Section 2).

Clearly, the complexity of the problem depends also on the presentation of the
function $f$: we will mostly assume that $f$ is presented by a {\em comparison
oracle} that, queried on $y,z\in\Z^d$, asserts whether $f(y)\leq f(z)$. This
is a broad presentation that reveals little information on the function, making
the problem harder to solve. In particular, if $d$ is variable, then already
for binary weights and maximizing a convex $f$, an exponential number of oracle
queries is needed (Proposition \ref{NPHard} part 3, Section 2).

In spite of these difficulties, we are able to provide the following efficient
algorithms for the problem: in the statements below, {\em oracle-time}
refers to the running time plus the number of oracle queries.

Our first theorem provides an efficient
algorithm for maximizing convex functions.
\bt{Convex}
For any fixed $d$, there is an algorithm that, given any positive integer $n$,
any integer weights $w^1,\dots,w^d$, and any convex function
$f:\R^d\longrightarrow\R$ presented by comparison oracle, solves the maximum
nonlinear bipartite matching problem in oracle-time which is polynomial
in $n$ and $\max|w^k_{i,j}|$.
\et

A second theorem provides an efficient
randomized algorithm for any function.
\bt{Random}
For any fixed $d$, there is a randomized algorithm that, given any positive
integer $n$, any integer weights $w^1,\dots,w^d$, and any function
$f:\R^d\longrightarrow\R$
presented by comparison oracle, solves the nonlinear bipartite matching
problem in oracle-time which is polynomial in $n$ and $\max|w^k_{i,j}|$.
\et

We also consider the minimum and maximum nonlinear bipartite matching
problems where the function $f$ is the $l_p$ norm
$\|\cdot\|_p:\R^d\longrightarrow\R$ given by
$\|y\|_p=(\sum_{k=1}^d |y_k|^p)^{1\over p}$ for $1\leq p<\infty$ and
$\|y\|_{\infty}=\max_{k=1}^d |y_k|$. For $l_p$ norm minimization we give
an algorithm which is polynomial in $n$ and $\max |w^k_{i,j}|$ and determines
a $d$-approximative solution for any $p$ and a more accurate,
$\sqrt d$-approximative solution, for the case of the Euclidean norm $p=2$
(Theorem \ref{MinNorm}). For $l_p$ norm maximization we give an algorithm
which is polynomial even in the bit size of the weights $w^k_{i,j}$ and
even if $d$ is variable, and determines a $d^{1\over p}$-approximative
solution for any $p$ (Theorem \ref{MaxNorm}).

The article proceeds as follows. In Section 2 we discuss various variants and
relatives of the problem, survey what is known in the literature about their
complexity, and demonstrate the intractability of the problem under various
conditions. In Section 3 we discuss convex maximization and prove Theorem
\ref{Convex}. In Section 4 we discuss approximative norm minimization
and maximization and prove Theorems \ref{MinNorm} and \ref{MaxNorm}. Finally,
in Section 5 we discuss randomized optimization and prove Theorem \ref{Random}.

\section{Variants and Intractability}

We now discuss various variants and relatives of nonlinear bipartite
matching, survey what is known (and unknown) about their complexity,
and demonstrate its intractability under various conditions.

First, consider the following related decision problem, asking
for the existence of a perfect matching attaining specified values
under each of the given weight (linear objective) functions $w^1,\dots,w^d$.

\vskip.2cm\noindent
{\bf Specified multiobjective bipartite matching}. Given $d,n$,
weight functions $w^1,\dots,w^d:N\longrightarrow\Z$, and integers
$u_1,\dots,u_d$, decide if there is a perfect matching
$M\subset N$ satisfying $w^k(M)=u_k$ for all $k$.

\vskip.2cm\noindent
Chandrasekaran et. al. considered the problem with a single objective $w$
(fixed $d=1$) and have shown that already this special case is
NP-complete \cite{CKM}. This raises the question about the complexity
in terms of the unary size $\max|w^k_{i,j}|$ of the weights. Indeed,
even the case of {\em binary} weights $w^k_{i,j}\in\{0,1\}$ is not
yet understood: for $d=1$ it was posed as intriguing and mysterious
by Papadimitriou and Yanakakis \cite{Pap,PY}, and the solutions obtained
consequently (first by Karzanov \cite{Kar} and recently in \cite{YMS})
are rather sophisticated; for $d=2$, the complexity is long open;
and for variable $d$ it is NP-complete.

The following proposition summarizes the known intractability
facts about the specified multiobjective bipartite matching problem.
We include the short proof for completeness of the exposition.

\bp{NPC}
The specified multiobjective bipartite matching problem is NP-complete already
under the following restrictions: (1) fixed $d=1$ (single weight function);
(2) binary weights $w^k_{i,j}\in\{0,1\}$.
\ep
\boproof
(1) is from \cite{CKM} by reduction from {\em subset sum}: given integers
$a_0,a_1,\dots,a_m$, define $d:=1$, $n:=2m$, single weight
$w\in\Z^{n\times n}$ by $w_{i,j}:=a_i$ for $1\leq i,j\leq m$ and $w_{i,j}:=0$
otherwise, and $u:=a_0$. Then there is a perfect matching $M$ with $w(M)=u$ if
and only if there is an $I\subseteq\{1,\dots,m\}$ with $a_0=\sum_{i\in I}a_i$.
\break
(2) is by reduction from {\em 3-dimensional matching}: given binary
${n\times n\times n}$ array $x$, put $d:=n$, $w^k_{i,j}:=x_{i,j,k}$, and
$u_k:=1$ for all $i,j,k$. Then there is a perfect matching $M$ with
$w^k(M)=u_k=1$ for all $k$ if and only if there is a binary array $y\leq x$ with
$\sum_{i,j}y_{i,j,k}=\sum_{i,k}y_{i,j,k}=\sum_{j,k}y_{i,j,k}=1$ for all $i,j,k$.
\eoproof

A further specialization of the case of binary weights $w^k_{i,j}\in\{0,1\}$
arises when the $w^k$ have pairwise disjoint supports. This can be formulated
as the following particularly appealing ``colorful" problem.

\vskip.2cm\noindent
{\bf Colorful bipartite matching}.
Given any bipartite graph $G$ with $d$-colored edge set
$E=\biguplus_{k=1}^d E_k$ and $u_1,\dots,u_d$, decide if there is a perfect
matching $M\subseteq E$ containing $u_k$ edges of color $k$ for $k=1,\dots,d$.

\vskip.2cm\noindent
This problem is a special case of specified multiobjective bipartite
matching with binary weights. To see this, note that we may assume
$G$ has the same number $n$ of vertices on each side,
making it a subgraph of $K_{n,n}$ with $E\subseteq N$,
and $\sum_{k=1}^d u_k=n$, else $G$ has no colorful perfect matching;
now, letting $w^k\in\{0,1\}^{n\times n}$ be the indicator of $E_k$
for all $k$, we have that $M\subset N$ is a perfect matching
of $K_{n,n}$ with $w^k(M)=u_k$ for all $k$ if and only if $M$
is a perfect matching of $G$ with $|M\cap E_k|=u_k$ for all $k$.

For $d=2$ (two colors), this problem is sometimes referred to in the
literature as the {\em exact matching problem}: for $G=K_{n,n}$ it is
polynomial time decidable \cite{Kar,YMS}; for arbitrary bipartite
graph $G$ there is a randomized algorithm \cite{MVV} but the
deterministic complexity is a longstanding open problem.

\vskip.3cm
Returning to nonlinear bipartite matching, the next proposition describes
its intractability under various conditions. By saying that an optimization
problem (rather than a decision problem) is NP-hard, we mean, as usual,
that there can be no polynomial time algorithm for solving it unless P=NP.

\bp{NPHard}
The following hold for the nonlinear bipartite matching problem
with data $d,n$, weights $w^1,\dots,w^d\in\Z^{n\times n}$, and function
$f:\R^d\longrightarrow\R$ presented explicitly or by a comparison oracle:
\begin{enumerate}
\item
For fixed $d=1$ (single weight function) and minimizing
the simple convex function $f_u:\R\longrightarrow\R$ defined by
$f_u(y):=(y-u)^2$ with $u$ an integer input parameter,
the problem is already NP-hard.
\item
For binary weights $w^k_{i,j}\in\{0,1\}$ and minimizing
the convex function $f_u:\R^d\longrightarrow\R$ defined by
$f_u(y):=\sum_{k=1}^d(y_k-u_k)^2$ with $u=(u_1,\dots,u_d)$
an integer vector, the problem is already NP-hard.
\item
For binary weights $w^k_{i,j}$ and maximizing a convex $f$
presented by comparison oracle, exponentially many oracle queries are
needed and hence the problem is not solvable in polynomial oracle-time.
\end{enumerate}
\ep

\boproof

\noindent
1.
Given weight $w$ and integer $u$, there is a perfect matching $M$ with
$w(M)=u$ if and only if the minimum value $f_u(w(M))=(w(M)-u)^2$ of a perfect
matching $M$ under $f_u$ is $0$. So even computing the optimal objective
function {\em value} enables deciding the NP-complete problem (1) of
Proposition \ref{NPC}.

\noindent
2.
Analogously to the proof of part 1 above: given binary weights $w^1,\dots,w^d$
and vector $u=(u_1,\dots,u_d)$, there is a perfect matching $M$ with
$w^k(M)=u_k$ for all $k$ if and only if the minimum objective value
$f_u(w^1(M),\dots,w^d(M))=\sum_{k=1}^d(w^k(M)-u_k)^2$ of a perfect matching
$M$ under $f_u$ is $0$. So even computing the optimal objective value
enables deciding the NP-complete problem (2) of Proposition \ref{NPC}.

\noindent
3.
Let $d=n^2$, define binary weights $w^{r,s}\in\Z^{n\times n}$ for
$1\leq r,s \leq n$, with $w^{r,s}_{i,j}:=1$ if $(i,j)=(r,s)$ and
$w^{r,s}_{i,j}=0$ otherwise, and let
$f:\R^d\cong\R^{n\times n}\longrightarrow\R$ be any function.
Then for any matrix $x\in\R^{n\times n}$ we have
$f(w^{1,1} x,\dots, w^{n,n} x)=f(x_{1,1},\dots, x_{n,n})=f(x)$.
Since the permutation matrices (which correspond to perfect matchings)
are convexly independent, any assignment of values to the $n!$
permutation matrices can be extended to a convex function $f$
on $\R^{n\times n}$. Thus, to find the permutation matrix
maximizing $f$, the oracle presenting $f$ must
be queried on all $n!$ permutation matrices.
\eoproof

\vskip.3cm
More generally, the {\em nonlinear combinatorial optimization problem} is
the following: given positive integers $d,n$, a family $\F$ of subsets of a
ground set $\{1,\dots,n\}$, integer weights $w^1,\dots,w^d$ on $\{1,\dots,n\}$,
and an arbitrary function $f:\R^d\longrightarrow\R$, find $F\in\F$ maximizing
(or minimizing) $f(w^1(F),\dots,w^d(F))$.
In \cite{OR}, the maximization problem with $f$ convex and $d$ fixed was
studied. It was shown that if the number of edge-directions of the polytopes
$P^{\F}:=\conv\{{\bf 1}_F:F\in\F\}$ (where ${\bf 1}_F\in\{0,1\}^n$
denotes the indicator of $F$) is polynomial in $n$ for a class of families
presented by membership oracles then the problem over families $\F$ in that
class can be solved in strongly polynomial oracle-time.
This unified and extended earlier results of \cite{Onn,OS} and yielded
polynomial time algorithms for convex maximization for various problems
including vector partitioning, matroids, and transportation problems with
fixed number of suppliers. However, for bipartite matching, which is the
combinatorial optimization problem over the family $\F\subset 2^N$ of perfect
matchings in $K_{n,n}$, the underlying polytope is the Birkhoff polytope
$$P^{\F}\ =\ \Pi^n\ :=\ \{x\in\R_+^{n\times n}:
\ \sum_i x_{i,j}=1, \ \sum_j x_{i,j}=1\}$$
which, as we next show, has exponentially many edge-directions,
making the methods of \cite{Onn,OR,OS} fail.
\bp{EdgeDirections}
The Birkhoff polytope $\Pi^n$ has precisely
${1\over2}\sum_{k=2}^n {n\choose k}^2 k!(k-1)!\geq{1\over n}{n!\choose 2}$
edge-directions.
\ep
\boproof
Every edge-direction of $\Pi^n$ is a nonzero minimal-support matrix
$x\in\R^{n\times n}$ with zero row-sums and column-sums,
and hence (up to scalar multiplication) is the matrix $x_C$
of some circuit $C\subset N$ of $K_{n,n}$, having values $\pm 1$
alternating along the edges of the circuit and $0$ elsewhere
(see e.g. \cite{OR}). We claim that each such circuit matrix $x_C$
is an edge-direction. To see this, let $C=C^+\biguplus C^-$ be the
partition of alternating edges of $C$ and let $D$ be a matching in
$K_{n,n}$ which perfectly matches all vertices not in $C$.
Let $x^+$ and $x^-$ be the permutation matrices which are
the indicators of the perfect matchings
$M^+:=C^+\cup D$ and $M^-:=C^-\cup D$ of
$K_{n,n}$. Define a binary weight matrix
$w$ as the indicator of $C\cup D$. Then $wx^+=wx^-=n$ whereas
$wx<n$ for any other permutation matrix $x$. Thus, $wx$ attains its
maximum over $\Pi^n$ precisely at the two vertices
$x^+$ and $x^-$ and hence $[x^+,x^-]$ is an edge and the difference
$x_C=x^+-x^-$ is an edge-direction. Now, for each $k\geq 2$, the
number of $2k$-circuits of $K_{n,n}$ is known and easily seen to be
${1\over2}{n\choose k}^2 k!(k-1)!$ (see \cite{OR}) and hence
the Proposition follows.
\eoproof

Proposition \ref{EdgeDirections} shows that, while the methods of
\cite{Onn,OR,OS} do apply for transportation problems with fixed number
of suppliers, they fail for bipartite matching which is the simplest
possible transportation problem - albeit, with variable numbers
of suppliers and consumers - and do not lead to a polynomial time
algorithm even for maximizing a convex $f$ with $d$ fixed.
This state of things, along with the easy solvability of the standard linear
bipartite matching problem, make the nonlinear problem for
bipartite matching particulary intriguing, and is part of our motivation
in raising and studying it herein.

\section{Deterministic Convex Maximization}

In this section we discuss the maximum nonlinear bipartite matching problem for
convex functions $f:\R^d\longrightarrow\R$ presented by comparison oracles.
We start with some definitions.
Here we will be working with matrices rather than graphs and matchings,
so the weights are now integer matrices $w^1,\dots,w^d\in\Z^{n\times n}$,
and the solutions are permutation matrices, which are well known to
be precisely the vertices of the Birkhoff polytope of bistochastic matrices
(with $\R_+$ the nonnegative reals),
$$\Pi^n\ =\ \{x\in\R_+^{n\times n}:
\ \sum_i x_{i,j}=1,\ \sum_j x_{i,j}=1\}\ .$$
Given weights $w^1,\dots,w^d$, define a projection
$w:\R^{n\times n}\longrightarrow\R^d$ mapping matrices $x$ to vectors
$w\cdot x$,
$$w\cdot x \quad:=\quad (w^1 x,\dots,w^d x)\quad=\quad
(\sum_{i,j}w^1_{i,j} x_{i,j},\dots,\sum_{i,j}w^d_{i,j} x_{i,j})\quad .$$
Define the {\em multiobjective polytope} (corresponding to
$w^1,\dots,w^d$) to be the projection of $\Pi^n$ under $w$,
$$\Pi^n_w\quad :=\quad
\{w\cdot x=(w^1 x,\dots,w^d x)\ :\ x\in\Pi^n\}\quad\subset\quad\R^d\quad .$$
Finally, define the {\em fiber} of any point $y\in\R^d$ to be the polytope
$\Pi^n\cap w^{-1}(y)$ consisting of those matrices in the Birkhoff polytope
that are projected by $w$ onto $y$. Thus, a point $y$ is in $\Pi_w^n$ if
and only if its fiber is nonempty; the following lemma asserts that
these equivalent conditions can be decided efficiently.

\bl{Fiber}
There is a polynomial time algorithm that, given $d,n,w^1,\dots,w^d$,
and integer $y\in\Z^d$, either asserts $y\not\in\Pi_w^n$ and
$\Pi^n\cap w^{-1}(y)=\emptyset$ or asserts $y\in\Pi_w^n$ and returns
a vertex $x$ of $\Pi^n\cap w^{-1}(y)$.
\el
\boproof
The fiber of any $y=(y_1,\dots,y_d)$ is the polytope
given by the following inequality description,
\begin{eqnarray*}
\Pi^n\cap w^{-1}(y) & = & \{x\in\Pi^n\ :\ (w^1 x,\dots,w^d x)
=(y_1,\dots,y_d)\, \} \\ & = & \{x\in\R_+^{n\times n}:
\ \sum_i x_{i,j}=1,\ \sum_j x_{i,j}=1,\ w^k x=y_k\}\ \ ,
\end{eqnarray*}
so linear programming allows to efficiently compute
a vertex of the fiber or assert that it is empty.
\eoproof

This lemma implies in turn that the multiobjective
polytope $\Pi_w^n$ can be constructed efficiently.

\bl{MultiobjectivePolytope}
For any fixed $d$, there is an algorithm that, given $n$ and
$w^1,\dots,w^d\in\Z^{n\times n}$, computes the vertex set
$\vert(\Pi_w^n)$ of the multiobjective polytope $\Pi_w^n$
in time which is polynomial in $n$ and $\max|w^k_{i,j}|$.
\el
\boproof
Let $u:=\max|w^k_{i,j}|$. Then for any permutation matrix $x$
and its projection $y=w\cdot x$ we have $|y_k|=|w^k x|\leq nu$ and
therefore $y$ lies in the grid $\{0,\pm 1,\dots,\pm nu\}^d$.
Since each vertex $y$ of $\Pi_w^n$ is the projection $y=w\cdot x$
of some vertex $x$ of $\Pi^n$, which is a permutation matrix,
we have $\vert(\Pi_w^n)\subseteq \{0,\pm 1,\dots,\pm nu\}^d$.
For each of the $(2nu+1)^d$ grid points $y\in\{0,\pm 1,\dots,\pm nu\}^d$,
apply the algorithm of Lemma \ref{Fiber} to check if $y\in\Pi_w^n$,
and obtain $Y:=\{0,\pm 1,\dots,\pm nu\}^d\cap\Pi_w^n$.

We then have that $\vert(\Pi_w^n)\subseteq Y\subseteq \Pi_w^n$ and
therefore the multiobjective polytope $\Pi_w^n$ is the convex hull of $Y$.
Since convex hulls can be computed in polynomial time for any fixed
dimension $d$, we can efficiently construct $\Pi_w^n$, that is,
determine all its vertices (and more generally all its faces).
\eoproof

Lemma \ref{Fiber} shows that for any $y\in\Z^d$
it is possible to check efficiently if $y$ is the projection
$y=w\cdot x$ of some bistochastic matrix $x\in\Pi^n$, and to
find such an $x$ if one exists. We will need also to
consider the integer analog of this problem: given $y\in\Z^d$,
is $y$ the projection $y=w\cdot x$ of some permutation matrix
$x\in\vert(\Pi^n)$, and if it is, can we find one such $x$
efficiently? but this problem is {\em precisely} the
specified multiobjective bipartite matching problem: there is
a permutation matrix $x$ with $y=w\cdot x$ if and only if there
is a perfect matching $M$ with $w^k(M)=y_k$ for $k=1,\dots,d$.
Unfortunately, as explained in Section 2, the complexity of this
problem is open even for fixed $d=2$. The difficulty is that the
fiber $\Pi^n\cap w^{-1}(y)$ of $y$ is not necessarily an integer
polytope and it may have some fractional (bistochastic) matrices
and some integer (permutation) matrices $x$ as its vertices.

Fortunately, as the next lemma shows, the fibers
of {\em vertices} of $\Pi_w^n$ are better behaved.
\bl{IntegerFiber} Let $y$ be any vertex of the multiobjective
polytope $\Pi_w^n$. Then the fiber $\Pi^n\cap w^{-1}(y)$ of $y$ is a
nonempty integer polytope all of whose vertices are permutation matrices.
Thus, the polynomial time algorithm of Lemma \ref{Fiber} applied to
$y\in\vert(\Pi_w^n)$ returns a permutation matrix $x$ satisfying $y=w\cdot x$.
\el
\boproof
It is well known and easy to see that if $Q$ is the image of a
polytope $P$ under an affine map $a$, then the preimage
$P\cap a^{-1}(F)=\{x\in P\,:\,a(x)\in F\}$ of any face $F$
of $Q$ is a face of $P$. Thus, if $y$ is a vertex of
$\Pi_w^n$ then its fiber $\Pi^n\cap w^{-1}(y)$,
which is the preimage under the map $w$ of the face $\{y\}$ of $\Pi_w^n$,
is a face of $\Pi_n$. Therefore, the vertices of the nonempty fiber of $y$,
one of which will be returned by the algorithm of Lemma \ref{Fiber}, are
precisely the vertices of $\Pi^n$ which are contained in that fiber.
\eoproof

We can now prove our first theorem, providing
an efficient algorithm for convex maximization.

\vskip.2cm\noindent{\bf Theorem \ref{Convex}}
For any fixed $d$, there is an algorithm that, given any positive integer $n$,
any integer weights $w^1,\dots,w^d$, and any convex function
$f:\R^d\longrightarrow\R$ presented by comparison oracle, solves the maximum
nonlinear bipartite matching problem in oracle-time which is polynomial
in $n$ and $\max|w^k_{i,j}|$.

\vskip.3cm
\boproof
Since $f$ is convex on $\R^d$ and $f(w^1(\cdot),\dots, w^d(\cdot))$
is convex on $\R^{n\times n}$, and the maximum of a convex function
over a polytope is attained at a vertex of the polytope,
we have the following equality,
\begin{eqnarray*}
\max\, \{f(w^1 x,\dots, w^d x):\ x\in\vert(\Pi^n)\}
& = & \max\, \{f(w^1 x,\dots, w^d x):\ x\in\Pi^n\} \\
& = & \max\, \{f(y):\ y\in\Pi_w^n\}
\ \ = \ \ \max\, \{f(y):\ y\in\vert(\Pi_w^n)\}\ \ .
\end{eqnarray*}
Apply the algorithm of Lemma \ref{MultiobjectivePolytope} and compute
$\vert(\Pi_w^n)$. By repeatedly querying the comparison oracle of $f$,
identify a vertex $y^*\in \vert(\Pi_w^n)$ attaining maximum value $f(y)$.
Now apply the algorithm of Lemma~\ref{Fiber} to $y^*$ and, as guaranteed by
Lemma \ref{IntegerFiber}, obtain a permutation matrix $x^*$ in the
fiber of $y^*$, so that $y^*=w\cdot x^*=(w^1 x^*,\dots,w^d x^*)$ and
$f(w^1 x^*,\dots,w^d x^*)=f(y^*)$. Since $y^*$ attains the maximum on the
right-hand side of the equation above, $x^*$ attains the maximum on the
left-hand side. Thus, the perfect matching of
$K_{n,n}$ corresponding to the permutation matrix $x^*$ is optimal.
\eoproof

The most time consuming part of the algorithm underlying Theorem \ref{Convex}
is the repeated use of linear programming for testing fibers of points in
the grid $\{0,\pm 1,\dots,\pm nu\}^d$ to construct $\vert(\Pi_w^n)$. There
are various ways of improving the algorithm in practice, but they do not
seem to improve the worst case complexity. We now describe such a variant
of the algorithm which will usually be much faster since it will
typically test the fibers of some but not all points in the grid.

\vskip.2cm\noindent{\bf A variant of the convex maximization algorithm.}
\begin{enumerate}
\item
Find the smallest grid containing $\vert(\Pi_w^n)$ by solving,
for $k=1,\dots,d$, the two linear programs
$$s_k:=\min \{w^k x:\sum_i x_{i,j}=\sum_j x_{i,j}=1,\ x\geq 0\}\,,
\ t_k:=\max \{w^k x:\sum_i x_{i,j}=\sum_j x_{i,j}=1,\ x\geq 0\}\,;$$
then $\vert(\Pi_w^n)$ is contained in the grid
$Z:=\{y\in\Z^d\,:\ s_k\leq y_k\leq t_k,\ k=1,\dots,d\}$.
\item
By repeatedly querying the comparison oracle of $f$, order the grid points
by nonincreasing value under $f$ and label them $y^1,\dots,y^{|Z|}$,
so that $Z=\{y^1,\dots,y^{|Z|}\}$ and $f(y^1)\geq\dots\geq f(y^{|Z|})$.
\item
Apply the algorithm of Lemma \ref{Fiber} to test the fiber
of each $y^i$ in order, until the first $y^k$
for which the vertex $x^*$ of its fiber $\Pi^n\cap w^{-1}(y^k)$
returned by the algorithm is a permutation matrix.
\item
Output the perfect matching of $K_{n,n}$
corresponding to the permutation matrix $x^*$.
\end{enumerate}

\noindent
We claim that $x^*$ is an optimal solution to the maximum convex
bipartite matching problem. Indeed, note that
$f^*:=\max\{f(y): y\in\Pi_w^n\}=\max\{f(y):y\in\vert(\Pi_w^n)\}$
equals the optimal objective function value
$\max\{f(w^1 x,\dots, w^d x):x\in\vert(\Pi^n)\}$
(see proof of Theorem \ref{Convex}); let $y^m\in\vert(\Pi_w^n)\subseteq Z$
be a vertex achieving that maximum value $f(y^m)=f^*$;
by Lemma \ref{IntegerFiber}, the algorithm of Lemma \ref{Fiber} applied to
$y^m\in\vert(\Pi_w^n)$ returns a permutation matrix and so $k\leq m$;
this implies $f^*\geq f(y^k)\geq f(y^m)=f^*$ and hence
$f^*=f(y^k)=f(w^1 x^*,\dots, w^d x^*)$;
therefore $x^*$ achieves the optimal objective function value.

\vskip0.5cm
We end this section with an example of a maximum convex bipartite
matching problem, demonstrating all notions and algorithms discussed
above, some of which will be also used in later sections.

\be{Example}
Consider the maximum convex bipartite matching problem with
the following data:
$$d=2,\quad n=4,
\quad w^1=\left(
\begin{array}{cccc}
  1 & 0 & 0 & 0 \\
  1 & 0 & 0 & 1 \\
  1 & 1 & 0 & 0 \\
  0 & 0 & 0 & 1 \\
\end{array}
\right),
\quad w^2=\left(
\begin{array}{cccc}
  1 & 1 & 0 & 1 \\
  0 & 0 & 1 & 1 \\
  0 & 0 & 1 & 0 \\
  1 & 1 & 0 & 0 \\
\end{array}
\right),
\quad f(y)=y_1^2+y_2^2
\quad .$$
By solving the linear programs minimizing and maximizing $w^k x$
over $\Pi^4$ for $k=1,2$ (step 1 of the algorithm above) we get $s_1=s_2=0$,
$t_1=3$, $t_2=4$, and so the smallest grid containing $\vert(\Pi_w^4)$
is $Z:=\{y\in\Z^2:0\leq y_1\leq 3,\ 0\leq y_2\leq 4\}
\subsetneq\{0,\pm 1,\dots,\pm 4\}^2$ which contains $20$ points.
Figure~\ref{Figure} below depicts this grid and indicates the objective
function value $f(y)=y_1^2+y_2^2$ of each grid point. Ordering the points
by decreasing value under $f$ (step 2 above) we get
$y^1=(3,4)\,,\ y^2=(2,4)\,,\ \dots\,,\ y^{20}=(0,0)$.
Testing fibers of the $y^i$ in order (step 3 above), the fiber of $y^1$ is
found empty whereas the fiber of $y^2$ is nonempty and is the first
for which the algorithm of Lemma \ref{Fiber} returns a permutation matrix
$$
x^*\quad =\quad\left(
\begin{array}{cccc}
  1 & 0 & 0 & 0 \\
  0 & 0 & 0 & 1 \\
  0 & 0 & 1 & 0 \\
  0 & 1 & 0 & 0 \\
\end{array}
\right)
\quad .$$
Thus, the corresponding matching $M^*:=\{(1,1),(2,4),(3,3),(4,2)\}$
of $K_{4,4}$ is an optimal solution.

Figure~\ref{Figure} also shows the multiobjective polytope $\Pi_w^4$ and
its vertex set $\vert(\Pi_w^4)$ computed by the algorithm of
Lemma \ref{MultiobjectivePolytope}: blue circles
are non-vertex grid points in $\Pi_w^4$
and green diamonds are vertices of $\Pi_w^4$. The optimal point
$y^2=(2,4)$ which is found either by the algorithm above or by the
algorithm of Theorem \ref{Convex} is the vertex of $\Pi_w^4$ attaining
maximum value under $f$ and is a red square. Of particular interest is
the blue point $y=(1,2)$ whose fiber $\Pi^n\cap w^{-1}(y)$
is a non-integer polytope with $30$ vertices (more than the $24$
of the Birkhoff polytope upstairs!), all of which are fractional, such as
$$
\left(\begin{array}{cccc}
 0   &   0   &   0   &   1   \\
 0   &   0.5 &   0.5 &   0   \\
 1   &   0   &   0   &   0   \\
 0   &   0.5 &   0.5 &   0   \\
\end{array}\right)
,
\left(\begin{array}{cccc}
 0   &   0.5 &   0   &   0.5 \\
 0.5 &   0   &   0.5 &   0   \\
 0   &   0.5 &   0   &   0.5 \\
 0.5 &   0   &   0.5 &   0   \\
\end{array}\right)
,
\left(\begin{array}{cccc}
 0   &   0.25    &   0.25    &   0.5 \\
 0.25    &   0.75    &   0   &   0   \\
 0.25    &   0   &   0.75    &   0   \\
 0.5 &   0   &   0   &   0.5         \\
\end{array}\right)
,
\left(\begin{array}{cccc}
 0   &   0.2 &   0.4 &   0.4 \\
 0   &   0.8 &   0   &   0.2 \\
 0.4 &   0   &   0.6 &   0   \\
 0.6 &   0   &   0   &   0.4 \\
\end{array}\right)
,
$$
indicating the difficulty of the specified multiobjective
and colorful bipartite matching problems.
\begin{figure}[!ht] \centering
\includegraphics[width=0.5\textwidth]{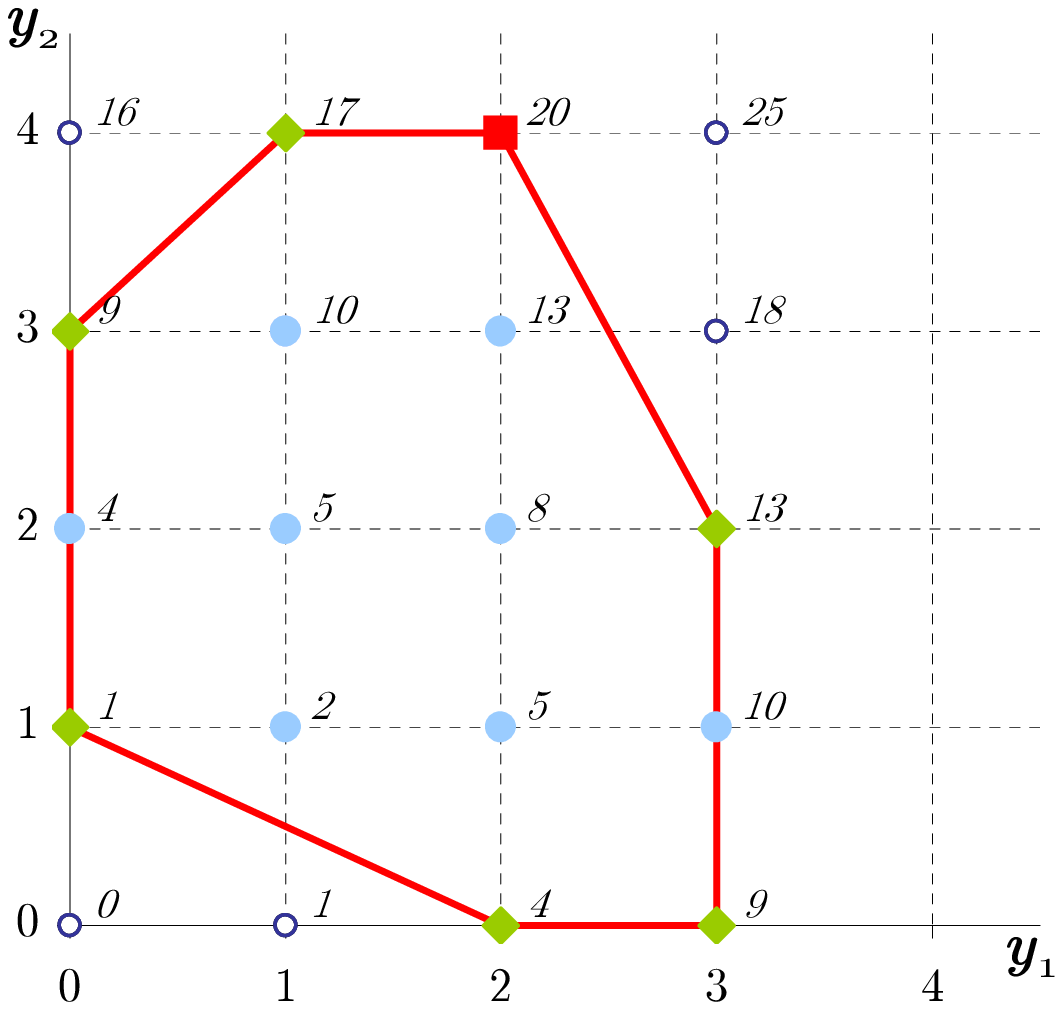}
\caption{Example \ref{Example} and the polytope $\Pi_w^4$}\label{Figure}
\end{figure}
\ee

\section{Approximative Norm Optimization}

Consider any discrete optimization problem with a finite set $S$ of feasible
solutions and nonnegative objective function $g:S\longrightarrow\R_+$
to be minimized or maximized and let $s^*\in S$ be any optimal solution.
Then an {\em $r$-approximative solution} is any feasible solution
$s\in S$ satisfying ${1\over r}g(s^*)\leq g(s)\leq r g(s^*)$.

In this section, building on the tools and results of Section 3, we provide
approximative algorithms for the minimum and maximum nonlinear
bipartite matching problems where the function $f$
is the $l_p$ norm $\|\cdot\|_p:\R^d\longrightarrow\R$ given by
$\|y\|_p=(\sum_{k=1}^d |y_k|^p)^{1\over p}$ for $1\leq p<\infty$ and
$\|y\|_{\infty}=\max_{k=1}^d |y_k|$. To keep our results general
and allow treatment of fractional and even nonrational $p$, we will
still assume that $f$ is presented by a comparison oracle. Of course,
for the most common values $p=1,2,\infty$ such an oracle is realizable
in time polynomial in the rest of the data; moreover, for any integer $p$,
by computing and comparing the integer valued $p$-th power $\|y\|^p_p$
of the norm instead of the norm itself, such an oracle is realizable
in time polynomial in the rest of the data and $\lceil \log p\rceil$.

\subsection{Minimization}

The following theorem provides an efficient
approximative algorithm for minimizing the $l_p$ norm.

\bt{MinNorm}
For any fixed $d$, there is an algorithm that, given any $n$, any
$1\leq p\leq\infty$, and any nonnegative integer weights $w^1,\dots,w^d$,
determines a $d$-approximative solution to the minimum nonlinear bipartite
matching problem with $f=\|\cdot\|_p$, in oracle-time which is polynomial
in $n$ and $\max w^k_{i,j}$. \break For $p=2$ (Euclidean norm), the
algorithm determines a more accurate, $\sqrt d$-approximative, solution.
\et
\boproof
The algorithm is the following: apply the algorithm of
Lemma \ref{MultiobjectivePolytope} and construct the vertex set
$\vert(\Pi_w^n)$ of the multiobjective polytope. Using the comparison
oracle of $f$ identify a vertex $\hat{y}\in \vert(\Pi_w^n)$
attaining minimum value $\|y\|_p$. Now apply the algorithm of
Lemma~\ref{Fiber} to $\hat{y}$ and, as guaranteed by
Lemma \ref{IntegerFiber}, obtain a permutation matrix $\hat{x}$
in the fiber of $\hat{y}$, so that $\hat{y}=w\cdot\hat{x}$.
Output the perfect matching of $K_{n,n}$ corresponding to
the permutation matrix $\hat{x}$.

We now show that this provides the claimed approximation. Let $x^*$ be the
permutation matrix corresponding to an optimal perfect matching and let
$y^*:=w\cdot x^*$ be its projection. Let $y'$ be a point on the boundary
of $\Pi_w^n$ satisfying $y'\leq y^*$. By Carath\'eodory's theorem
(on the boundary) $y'$ is a convex combination $y'=\sum_{i=1}^r\lambda_i y^i$
of some $r\leq d$ vertices of $\Pi_w^n$ and hence
$\lambda_t=\max \lambda_i\geq{1\over r}\geq {1\over d}$ for some $t$.

Since the weights $w^k$ are nonnegative we find that
so are $y'$ and the $y^i$ and hence we obtain
\begin{eqnarray*}
f(w^1\hat{x},\dots,w^d\hat{x})& = & \|\hat{y}\|_p \quad \leq\quad \|y^t\|_p
\quad\leq\quad d\lambda_t\cdot\|y^t\|_p\quad=\quad d\cdot\|\lambda_t y^t\|_p \\
& \leq & d\cdot\left\|\sum_{i=1}^r\lambda_i y^i\right\|_p\quad =\quad
d\cdot\|y'\|_p \quad \leq\quad d\cdot\|y^*\|_p
\quad =\quad d\cdot f(w^1x^*,\dots,w^dx^*)\quad.
\end{eqnarray*}
This proves that $\hat{x}$ provides a $d$-approximative solution for
any $1\leq p\leq\infty$. Now consider the case of Euclidean norm $p=2$.
By Cauchy-Schwartz,
$1=(\sum_{i=1}^r 1\cdot\lambda_i)^2\leq \sum_{i=1}^r 1^2
\sum_{i=1}^r\lambda_i^2=r\sum\lambda_i^2\leq d\sum\lambda_i^2$.
Find $s$ with $\|y^s\|_p=\min\|y^i\|_p$ and recall that the $y^i$
are nonnegative. We then have the inequality
\begin{eqnarray*}
f^2(w^1\hat{x},\dots,w^d\hat{x})& = & \|\hat{y}\|_2^2 \quad\leq\quad\|y^s\|_2^2
\quad\leq\quad (d\sum_{i=1}^r\lambda_i^2)\cdot\|y^s\|_2^2\quad\leq\quad
d\sum_{i=1}^r\lambda_i^2 \cdot \|y^i\|_2^2 \\
& \leq & d\cdot\left\|\sum_{i=1}^r\lambda_i y^i\right\|_2^2
\quad = \quad d\cdot\|y'\|_2^2 \quad \leq\quad d\cdot\|y^*\|_2^2
\quad =\quad d\cdot f^2(w^1x^*,\dots,w^dx^*)
\end{eqnarray*}
which proves that in this case, as claimed, $\hat{x}$
provides moreover a $\sqrt d$-approximative solution.
\eoproof

\subsection{Maximization}

The following theorem provides an approximative algorithm for maximizing
the $l_p$ norm, that runs in time which is polynomial even in the
bit size of the weights $w^k_{i,j}$ and even if $d$ is variable.

\bt{MaxNorm}
There is an algorithm that, given any $d$, any $n$, any $1\leq p\leq\infty$,
and any nonnegative integer weights $w^1,\dots,w^d$, determines a
$d^{1\over p}$-approximative solution to the maximum nonlinear bipartite
matching problem with $f=\|\cdot\|_p$, in oracle-time which is polynomial
in $d$ ,$n$, and $\max\lceil \log w^k_{i,j}\rceil$.
\et
\boproof
The algorithm is the following:
for $k=1,\dots,d$ solve the linear programming problem
$$\max\{w^k x\ :\ \sum_i x_{i,j}=1,\ \sum_j x_{i,j}=1,\ x\geq 0\}\ ,$$
obtain an optimal vertex $x^k$ of $\Pi^n$, and let $y^k:=w\cdot x^k$ be its
projection. Using the comparison oracle of $f$ find $r$ with
$\|y^r\|_p=\max_{k=1}^d\|y^k\|_p$.
Output the perfect matching of $K_{n,n}$ corresponding to $x^r$.

We now show that this provides the claimed approximation.
Let $s$ satisfy $\|y^s\|_\infty=\max_{k=1}^d\|y^k\|_\infty$.
First, we claim that any $y\in\Pi^n_w$ satisfies
$\|y\|_\infty\leq \|y^s\|_\infty$. To see this, choose any point
$x\in\Pi^n\cap w^{-1}(y)$ in the fiber of $y$ so that $y:=w\cdot x$,
let $t$ satisfy $y_t=\|y\|_\infty=\max_{k=1}^d y_k$, and recall
that the $w^k$ and hence the $y^k$ are all nonnegative.
Then, as claimed, we get
$$\|y\|_\infty \ \ = \ \ y_t \ \ =\ \ w^t x
\ \ \leq \ \ \max\{w^t x\,:\, x\in\Pi^n\} \ \ = \ \ w^t x^t \ \ = \ \ y^t_t
\ \ \leq \ \ \|y^t\|_\infty \ \ \leq\ \ \|y^s\|_\infty \ \ .$$

Let $x^*$ be an optimal permutation matrix
and let $y^*:=w\cdot x^*$ be its projection.
Consider first the case $p=\infty$. Then
$\|y^r\|_\infty=\max_{k=1}^d\|y^k\|_\infty=\|y^s\|_\infty$
and hence, by the claim just proved, we have
$$f(w^1 x^*,\dots,w^d x^*) \ \ =\ \ \|y^*\|_\infty\ \ \leq
\ \ \|y^s\|_\infty \ \ = \ \ \|y^r\|_\infty\ \ = \ \
f(w^1 x^r,\dots,w^d x^r) \ \ \leq \ \ f(w^1 x^*,\dots,w^d x^*)\ \ .$$
Therefore equality holds all along and $x^r$ provides an exact optimal
solution, or in other words, a $1$-approximative solution, agreeing with
the statement of the theorem with $d^{1\over\infty}=1$ for $p=\infty$.
Next, consider the case of any $1\leq p<\infty$.
Then we have the following inequality which completes the proof,
\begin{eqnarray*}
f^p(w^1 x^*,\dots,w^d x^*)& = & \|y^*\|_p^p
\quad = \quad \sum_{k=1}^d |y^*_k|^p \quad \leq \quad d\cdot \|y^*\|^p_\infty
\quad \leq \quad d\cdot \|y^s\|^p_\infty \\
& \leq & d\cdot \sum_{k=1}^d |y^s_k|^p
\quad = \quad d\cdot \|y^s\|_p^p \quad \leq \quad d\cdot \|y^r\|_p^p
\quad = \quad  d\cdot f^p(w^1x^s,\dots,w^dx^s)
\ .\ \square
\end{eqnarray*}

\section{Randomized Nonlinear Optimization}

In this section we provide a randomized algorithm for nonlinear
bipartite matching for any function $f:\R^d\longrightarrow\R$
presented by a comparison oracle. By this we mean an algorithm
that has access to a random bit generator, and on any input
outputs the optimal solution with probability at least half.

By adding to each $w^k_{i,j}$ a suitable positive integer $v$
and replacing the function $f$ by the function that maps each
$y\in\R^d$ to $f(y_1-nv,\dots,y_d-nv)$ if necessary, we may and will
assume without loss of generality throughout this section that
the given weights are nonnegative, $w^1,\dots,w^d\in\N^{n\times n}$.

Recall that $\vert(\Pi^n)$ is the set of $n\times n$ permutation matrices
and let $Y:=\{w\cdot x\,:\,x\in\vert(\Pi^n)\}$ be the set
of all projections $y=w\cdot x=(w^1 x,\dots,w^d x)\in\N^d$
of permutation matrices $x$.

In this section we will be working with polynomials with
integer coefficients in the $n^2+d$ variables $a_{i,j}$, $i,j=1,\dots,n$
and $b_k$, $k=1,\dots,d$.
Define an $n\times n$ matrix $A$ whose entries are monomials by
$$A_{i,j}\quad:=\quad
a_{i,j}\prod_{k=1}^d b_k^{w_{i,j}^k},\quad i,j=1,\dots,n\quad.$$
For each matrix $x\in\N^{n\times n}$ and vector $y\in\N^n$,
the corresponding monomial is
$$a^xb^y\quad:=\quad\prod_{i=1}^n\prod_{j=1}^n a_{i,j}^{x_{i,j}}
\prod_{k=1}^d b_k^{y_k}\quad .$$
For each permutation matrix $x$ let $\sign(x)=\pm$ denote the sign of the
corresponding permutation. Finally, for each $y\in\N^d$ define
the following polynomial in the variables $a=(a_{i,j})$ only, by
$$g_y(a)\quad:=\quad\sum
\left\{\sign(x)a^x\,:\,x\in\vert(\Pi^n),\ \ w\cdot x=y\right\}\quad.$$
We then have the following identity expanding the
determinant of $A$ in terms of the $g_y(a)$,
$$\det(A)\ \ =\ \ \sum_{x\in\vert(\Pi^n)}
\sign(x)\prod_{i,j}A_{i,j}^{x_{i,j}}
\ \ =\ \ \sum_{x\in\vert(\Pi^n)}\sign(x)a^xb^{w\cdot x}
\ \ =\ \ \sum_{y\in Y=w\cdot\vert(\Pi^n)}g_y(a)b^y\ \ .
$$

Next we consider integer substitutions to the variables $a_{i,j}$.
Under such substitutions, each $g_y(a)$ becomes an integer and
$\det(A)=\sum_{y\in Y}g_y(a)b^y$ becomes a polynomial in the variables
$b=(b_k)$ only. Given such a substitution,
let $\hat{Y}:=\{y\in Y\,:\,g_y(a)\neq 0\}$ be the {\em support}
of $\det(A)$, that is, the set of exponents of monomial $b^y$ appearing
with nonzero coefficient in $\det(A)$.

The next proposition concerns substitutions of independent identical
random variables uniformly distributed on the set of integers
$\{1,2,\dots,s\}$, under which $\hat{Y}$ becomes a random subset of $Y$.
\bp{Probability}
Suppose that independent identical random variables uniformly
distributed on the set $\{1,2,\dots,s\}$ are substituted for the $a_{i,j}$
and let $\hat{Y}=\{y\in Y\,:\,g_y(a)\neq 0\}$ be the random support of $\det(A)$.
Then, for every $y\in Y=\{w\cdot x\,:\,x\in\vert(\Pi^n)\}$,
the probability that $y\notin\hat{Y}$ is at most ${n\over s}$.
\ep
\boproof
Consider any $y\in Y$ and consider $g_y(a)$ as a polynomial in the variables
$a=(a_{i,j})$. Since $y=w\cdot x$ for some permutation matrix,
there is as least one term $\sign(x)a^x$ in $g_y(a)$.
Since distinct permutation matrices $x$ give distinct monomials $a^x$,
no cancellations occur among the terms $\sign(x)a^x$ in $g_y(a)$.
Thus, $g_y(a)$ is a nonzero polynomial of degree $n$.
The claim now follows from a lemma of Schwartz \cite{Sch}
stating that the substitution of independent identical
random variables uniformly distributed on $\{1,2,\dots,s\}$ into a nonzero
multivariate polynomial of degree $n$ is zero with probability
at most $n\over s$.
\eoproof

The next lemma shows that, given $a_{i,j}$, the support
$\hat Y$ of $\det(A)$ is polynomial time computable.

\bl{Interpolation}
For any fixed $d$, there is an algorithm that, given $n$,
$w^1,\dots,w^d\in\N^{n\times n}$, and substitutions $a_{i,j}\in\{1,2,\dots,s\}$,
computes $\hat{Y}=\{y\in Y\,:\,g_y(a)\neq 0\}$ in time polynomial in
$n$, $\max w^k_{i,j}$ and $\lceil \log s \rceil$.
\el
\boproof
For each $y$, let $g_y:=g_y(a)$ be the fixed integer obtained by
substituting the given integers $a_{i,j}$. Put $u= n\cdot \max w^k_{i,j}$
and $Z=\{0,1,\dots, u\}^d$. Then $\hat{Y}\subseteq Y\subseteq Z$
and hence $\det(A)=\sum_{y\in Z}g_yb^y$ is a polynomial in $d$ variables
$b=(b_k)$ involving at most $|Z|=(u+1)^d$ monomials. For $t=1,2,\dots,(u+1)^d$
consider the substitution $b_k:=t^{(u+1)^{k-1}}$, $k=1,\dots,d$.
Let $A(t)$ be the integer matrix obtained from $A$ by
this substitution along with the substitution of the given $a_{i,j}$.
Then each entry of $A(t)$ satisfies
$$A(t)_{i,j}\ \ =\ \ a_{i,j}\prod_{k=1}^d (t^{(u+1)^{k-1}})^{w_{i,j}^k}
\ \ \leq\ \ s\prod_{k=1}^d (((u+1)^d)^{(u+1)^{k-1}})^{u\over n}
\ \ \leq\ \ s(u+1)^{d(u+1)^{d+1}}
$$
and hence its bits size $1+\log A(t)_{i,j}=O(u^{d+1}\log(su))$ is
polynomially bounded in $n,\max w^k_{i,j},\lceil \log s \rceil$.
Therefore the integer number
$\det(A(t))$ can be computed in polynomial time by Gaussian elimination.
So we obtain the following system of $(u+1)^d$ equations in
$(u+1)^d$ variables $g_y$, $y\in Z=\{0,1,\dots, u\}^d$,
$$\det(A(t))\quad=\quad\sum_{y\in Z}g_y\prod_{k=1}^d b_k^{y_k}
\quad=\quad\sum_{y\in Z}t^{\sum_{k=1}^d y_k (u+1)^{k-1}}\cdot g_y
\,,\quad\quad t=1,2,\dots,(u+1)^d\quad.
$$
As $y=(y_1,\dots,y_d)$ runs through $Z$, the sum
$\sum_{k=1}^d y_k (u+1)^{k-1}$ attains precisely all $|Z|=(u+1)^d$
distinct values $0,1,\dots,(u+1)^d-1$. This implies that, under the
total order of the points $y$ in $Z$ by increasing value of
$\sum_{k=1}^d y_k (u+1)^{k-1}$, the vector of coefficients of the
$g_y$ in the equation corresponding to $t$ is precisely the point
$(t^0,t^1,\dots,t^{(u+1)^d-1})$ on the moment curve in
$\R^Z\simeq\R^{(u+1)^d}$. Therefore, the equations are linearly
independent and hence the system can be solved for the $g_y=g_y(a)$
and the desired support $\hat{Y}=\{y\in Y\,:\,g_y(a)\neq 0\}$
of $\det(A)$ can indeed be computed in polynomial time.
\eoproof

We are now in position to prove Theorem \ref{Random}.
By a {\em randomized algorithm} that solves the nonlinear bipartite matching
problem we mean an algorithm that has access to a random bit generator
and on any input to the problem outputs a perfect matching which
is optimal with probability at least a half. The running time
of the algorithm includes a count of the number of random bits used.
Note that by repeatedly applying such an algorithm several times and
picking the best perfect matching, the probability of failure can be
decreased at will; in particular, repeating it $n$ times decreases the
failure probability to as negligible a fraction as $1\over 2^n$ while
increasing the running time by a linear factor only.

\vskip.2cm\noindent{\bf Theorem \ref{Random}}
For any fixed $d$, there is a randomized algorithm that, given any positive
integer $n$, any integer weights $w^1,\dots,w^d$, and any function
$f:\R^d\longrightarrow\R$
presented by comparison oracle, solves the nonlinear bipartite matching
problem in oracle-time which is polynomial in $n$ and $\max|w^k_{i,j}|$.

\vskip.3cm
\boproof
As explained in the beginning of this section, we may and will assume
that the $w^k$ are nonnegative.
First we claim that, with probability at least $1-{1\over 2n}$,
we can compute the optimal {\em objective function value} of
the nonlinear bipartite matching problem.
To see this, note that the optimal value equals $\max\{f(y):y\in Y\}$
where $Y=\{w\cdot x\,:\,x\in\vert(\Pi^n)\}$ as before, and let
$y^*\in Y$ be a point attaining $f(y^*)=\max\{f(y):y\in Y\}$.
Now, using polynomially many random bits, draw independently and uniformly
distributed integers from $\{1,2,\dots,2n^2\}$ and substitute them for
the $a_{i,j}$. Next compute $\hat{Y}=\{y\in Y\,:\,g_y(a)\neq 0\}$ using
the algorithm underlying Lemma \ref{Interpolation} and determine
$\max\{f(y):y\in \hat{Y}\}$. By Proposition \ref{Probability},
with probability at least $1-{1\over 2n}$ we have $y^*\in\hat{Y}$
in which event $\max\{f(y):y\in \hat{Y}\}=\max\{f(y):y\in Y\}$
is indeed the optimal objective function value.

Next, suppose that $M\subset N=\{(i,j):1\leq i,j\leq n\}$ is any
(not necessarily perfect) matching of $K_{n,n}$. Then we can also
compute, with probability at least $1-{1\over 2n}$, the maximum objective
function value among perfect matchings $M\cup L$ containing $M$.
To see this, let $m:=n-|M|$ and consider the subgraph $G$
of $K_{n,n}$ induced by the vertices not matched under $M$.
Then $G$ is isomorphic to $K_{m,m}$ and we have a naturally
induced nonlinear bipartite matching problem on $G$,
where the new weight functions $\bar{w}^k$ are simply the restrictions
of the $w^k$ to the edges of $G$, and the new functional $\bar f$ on $\R^d$
is defined by $\bar{f}(y_1,\dots,y_d):=f(y_1+w^1(M),\dots,y_d+w^d(M))$.
Then the objective function value $f(w^1(M\cup L),\dots,w^d(M\cup L))$
of any perfect matching $M\cup L$ of $K_{n,n}$ in the original problem equals
the objective function value $\bar{f}(\bar{w}^1(L),\dots,\bar{w}^d(L))$
of the perfect matching $L$ of $G$ in the induced problem.
Since $\max\bar{w}^k_{i,j}\leq\max w^k_{i,j}$ and $m\leq n$ we can
compute, with probability at least $1-{1\over 2n}$, in time polynomial
in $n$ and $\max|w^k_{i,j}|$, the optimal objective function value of
a perfect matching of $G$ by the algorithm of the paragraph above applied
to $G$, where the randomized substitutions are taken from
$\{1,2,\dots,2mn\}$ (and not from $\{1,2,\dots,2m^2\}$ which would
give smaller probability of success). This value is the maximum objective
function value among perfect matchings $M\cup L$ containing $M$.

We claim that the following procedure constructs a perfect matching
$M$ of $K_{n,n}$ which is optimal with probability at least $1\over2$.
Start with $M:=\emptyset$ and $i:=1$. While $i\leq n$ iterate:
for each edge $(i,j)$ such that $j$ is not matched under $M$, use the
algorithm of the previous paragraph to obtain the maximum objective function
value of a perfect matching of $K_{n,n}$ containing $M\cup\{(i,j)\}$;
let $j_i$ be the smallest $j$ for which this value is maximal;
update $M:=M\cup\{(i,j_i)\}$; increment $i$ and repeat. Output $M$.

To prove the claim, let $M^*=\{(1,r_1),(2,r_2),\dots,(n,r_n)\}$ be the
{\em lexicographically first} optimal perfect matching of $K_{n,n}$, that is,
the one such that for any other optimal $M'=\{(1,s_1),(2,s_2),\dots,(n,s_n)\}$
there is an index $1\leq h<n$ such that $r_i=s_i$ for all $i<h$ and $r_h<s_h$.
For $i=1,\dots,n$ let $E_i$ be the random event that after the completion
of iteration $i$ of the above procedure we have $M=\{(1,r_1),\dots,(i,r_i)\}$.
We prove by induction on $i$ that $\Pr(E_i)\geq (1-{1\over 2n})^i$. For the
basis note that $E_1$ is the event that the randomized algorithm used during
the first iteration computes correctly the maximum objective function value
of a perfect matching containing $\{(1,r_1)\}$, having probability at least
$1-{1\over 2n}$. For the inductive step note that $\Pr(E_{i+1}|E_i)$ is the
probability that, given that after iteration $i$ we have
$M=\{(1,r_1),\dots,(i,r_i)\}$, the randomized algorithm used during iteration
$i+1$ computes correctly the maximum objective function value of a perfect
matching containing $\{(1,r_1),\dots,(i,r_i),(i+1,r_{i+1})\}$, which is again
at least $1-{1\over 2n}$; as $E_{i+1}\subseteq E_i$, the induction follows by
$$\Pr(E_{i+1})\ \ =\ \ \Pr(E_{i+1}\cap E_i)\ \ =\ \ \Pr(E_{i+1}|E_i)\Pr(E_i)
\ \ \geq\ \ (1-{1\over 2n})(1-{1\over 2n})^i\ \ =\ \ (1-{1\over 2n})^{i+1}\ \ .$$
Now, the probability that the perfect matching $M$ output by the procedure
above is optimal is no smaller than the probability that $M$ equals the
lexicographically first optimal perfect matching $M^*$ which is precisely
$\Pr(E_n)$ and hence at least $(1-{1\over 2n})^n\geq {1\over 2}$ as desired.
This completes the proof.
\eoproof

\noindent {\small Yael Berstein}\newline
\emph{Technion - Israel Institute of Technology, 32000 Haifa, Israel}\newline
\emph{email: yaelber{\small @}tx.technion.ac.il}

\vskip.5cm

\noindent {\small Shmuel Onn}\newline
\emph{Technion - Israel Institute of Technology, 32000 Haifa, Israel}\newline
\emph{email: onn{\small @}ie.technion.ac.il},
\ \ \emph{http://ie.technion.ac.il/{\small $\sim$onn}}

\end{document}